\input amstex
\documentstyle{amsppt}
\document

\magnification 1100

\def\gen{{\frak{g}}}

\def\sen{\frak{s}}

\def\len{{\frak l}}

\def\1b{{\bold 1}}
\def\eb{{\bold e}}
\def\fb{{\bold f}}

\def\hb{{\bold h}}

\def\X{{\roman X}}

\def\L{{\roman L}}

\def\SL{{\S\L}}
\def\S{{\roman S}}

\def\X{{\roman X}}

\def\Lie{\text{Lie}}

\def\Hom{\text{Hom}\,}

\def\res{\text{res}\,}

\def\Gr{\text{Gr}}

\def\CC{{\Bbb C}}

\def\GG{{\Bbb G}}

\def\PP{{\Bbb P}}

\def\ZZ{{\Bbb Z}}

\def\Cc{{\Cal C}}

\def\Ec{{\Cal E}}
\def\Fc{{\Cal F}}

\def\Ic{{\Cal I}}
\def\ICc{{\Ic\Cc}}

\def\Lc{{\Cal L}}

\def\Pc{{\Cal P}}

\def\and{{\quad\text{and}\quad}}

\def\qed{\hfill $\sqcap \hskip-6.5pt \sqcup$}        
\overfullrule=0pt                                    

\newdimen\Squaresize\Squaresize=14pt
\newdimen\Thickness\Thickness=0.5pt
\def\Square#1{\hbox{\vrule width\Thickness
	      \vbox to \Squaresize{\hrule height \Thickness\vss
	      \hbox to \Squaresize{\hss#1\hss}
	      \vss\hrule height\Thickness}
	      \unskip\vrule width \Thickness}
	      \kern-\Thickness}
\def\Vsquare#1{\vbox{\Square{$#1$}}\kern-\Thickness}

\title On the action of the dual group on the cohomology of
perverse sheaves on the affine grassmannian\endtitle
\rightheadtext{perverse sheaves on the affine grassmanian}
\author E. Vasserot\endauthor
\thanks
The author is partially supported by EEC grant
no. ERB FMRX-CT97-0100.\endthanks
\abstract  It was proved by Ginzburg and Mirkovic-Vilonen that the $G(O)$-equivariant
perverse sheaves on the affine grassmannian of a connected reductive group $G$
form a tensor category equivalent to the tensor category of finite dimensional
representations of the dual group $G^\vee$. The proof use the
Tannakian formalism. The purpose of this paper is to construct explicitely the
action of $G^\vee$ on the global cohomology of a perverse sheaf.
It would be interesting to find a $q$-analogue of this construction.
It would give the global counterpart to \cite{BG}.
\endabstract
\endtopmatter
\document
\vskip1cm

\head 1. Notations and reminder on affine Grassmanians\endhead
\subhead 1.1\endsubhead
Let $G$ be a connected reductive complex algebraic group.
Let $B, T$, be a Borel and a Cartan subgroup of $G$.
Let $U\subset B$ be the unipotent radical of $B$.
Let $B^-$ a Borel subgroup such that $B\cap B^-=T$.
Set $X_T=\Hom(T,\GG_m)$ and $X^\vee_T=\Hom(\GG_m,T)$
be the weight an the coweight lattice of $G$.
For simplicity we write $X=X_T$ and $X^\vee=X^\vee_T$.
Let $(\ ,\ )\,:\,X\times X^\vee\to\ZZ$ be the natural pairing.
Let $R$ be the set of roots, $R^\vee$ the set of coroots.
Let $R_\pm\subset R$, $R^\vee_\pm\subset R^\vee$, be the subsets
of positive and negative roots and coroots.
Let $X_+\subset X$, $X^\vee_+\subset X^\vee$, be the subsets
of dominant weights and coweights. Let $\rho_G\in X$ be half the sum
of all positive roots. 
If there is no ambiguity we simply write $\rho$ instead of 
$\rho_G$. Let $G^\vee$ and $Z(G)$ be the dual group and the center of $G$.
Let $\alpha_i,\alpha_i^\vee,$ $i\in I$, be the simple roots and the simple 
coroots, and let $\omega_i,\omega_i^\vee$ be the fundamental weights and 
coweights. For any root $\alpha\in R$ let $U_\alpha\subset G$ be the 
corresponding root subgroup. If $\alpha=\alpha_i$, $i\in I$, we simply
set $U_i=U_{\alpha_i}$ and $U^-_i=U_{-\alpha_i}$. 
Let $W$ be the Weyl group of $G$. For any
$i\in I$ let $s_i$ be the simple reflexion corresponding
to the simple root $\alpha_i$. 

\subhead 1.2\endsubhead
Let $K=\CC((t))$ be the field of Laurent formal series, and let
$O=\CC[[t]]$ be the subring of integers. 
Recall that $G(O)$ is a group scheme and that $G(K)$ is a group ind-scheme.
The quotient set
$\Gr^G=G(K)/G(O)$ is endowed with the structure of an ind-scheme.
We may write $\Gr$ instead of $\Gr^G$, hoping that it makes no confusion.
For any coweight $\lambda^\vee\in X^\vee$, let $t^{\lambda^\vee}\in T(K)$
be the image of $t$ by the group homomorphism 
$\lambda^\vee\,:\,\GG_m(K)\to T(K)$. If $\lambda^\vee$ is dominant,
set $e_{\lambda^\vee}=t^{\lambda^\vee}G(O)/G(O)\in\Gr.$
The $G(O)$-orbit $\Gr_{\lambda^\vee}=G(O)\cdot e_{\lambda^\vee}$
is connected and simply connected.
Let $\overline\Gr_{\lambda^\vee}$ be its Zariski closure.
Let $\Pc_G$ be the category of $G(O)$-equivariant perverse sheaves
on $\Gr$. For any $\lambda^\vee$ let $\ICc_{\lambda^\vee}$
be the intersection cohomology complex on $\Gr_{\lambda^\vee}$
with coefficients in $\CC$.
Consider the fiber product $G(K)\times_{G(O)}\Gr$.
It is the quotient of $G(K)\times\Gr$ by $G(O)$, where $u\in G(O)$
acts on $G(K)\times\Gr$ by $(g,x)\mapsto (gu^{-1},ux)$. 
The map 
$$\tilde p\,:\,G(K)\times_{G(O)}\Gr\to\Gr,\quad (g,x)\mapsto ge_0$$ 
is the locally trivial fibration with fiber $\Gr$
associated to the $G(O)$-bundle 
$$p\,:\,G(K)\to\Gr.$$ 
Thus $G(K)\times_{G(O)}\Gr$ is an ind-scheme :
it is the inductive limit of the subschemes 
$p^{-1}(\overline\Gr_{\lambda^\vee_1})\times_{G(O)}
\overline\Gr_{\lambda_2^\vee}.$
Consider also the map 
$$m\,:\,G(K)\times_{G(O)}\Gr\to\Gr,\quad (g,x)\mapsto gx.$$ 
For any $\lambda^\vee_1,\lambda_2^\vee\in X^\vee_+$
let $\ICc_{\lambda_1^\vee}\star\ICc_{\lambda_2^\vee}$ be the direct image
by $m$ of the intersection cohomology complex of the subvariety
$$p^{-1}(\overline\Gr_{\lambda^\vee_1})\times_{G(O)}
\overline\Gr_{\lambda_2^\vee}\subset G(K)\times_{G(O)}\Gr.$$
The complex $\ICc_{\lambda_1^\vee}\star\ICc_{\lambda_2^\vee}$ is
perverse (see \cite{MV}, and \cite{NP, Corollaire 9.7} for more details). 
It is known that the cohomology sheaves of the complex $\Ic\Cc_{\lambda^\vee}$
are pure by the argument similar to \cite{KT} (see also \cite{KL}).
It is also known that any object in $\Pc_G$ is a direct sum of 
complexes $\ICc_{\lambda^\vee}$ (see \cite{BD, Proposition 5.3.3.(i)}
for a proof). Thus we get a convolution product
$\star\,:\,\Pc_G\times\Pc_G\to\Pc_G$. It is the convolution product defined
by Mirkovic and Vilonen.

\subhead 1.3\endsubhead Let $P\subset G$ be a parabolic subgroup of $G$,
$N\subset P$ be the unipotent radical, $M=P/N$ be the Levi factor.
Let $M'=[M,M]$ be the semisimple part of $M$. Consider the diagram
$$\Gr^G\,{\buildrel\gamma\over\leftarrow}\,\Gr^P\,
{\buildrel\pi\over\to}\,\Gr^M,$$
where the maps $\gamma$ and $\pi$ are induced by the embedding $P\subset G$ and
the projection $P\to M$.
The fibers of $\pi$ are $N(K)$-orbits. 
Observe that $\Gr^M$ is not connected.
The connected components of $\Gr^M$ are labelled by characters of the center
of the dual group $M^\vee$.
Let $\Gr^{M,\theta^\vee}\subset\Gr^M$ be the component associated to 
$\theta^\vee\in X_{Z(M^\vee)}$.
By definition, $e_{\lambda^\vee}\in\Gr^{M,\theta^\vee}$
if and only if the restriction of $\lambda^\vee$ to $Z(M^\vee)$ coincides
with $\theta^\vee$.
The element $\rho-\rho_M$ belongs to $X^\vee_{Z(M^\vee)}$. Put
$$\Gr^{M,n}=\bigsqcup_{2(\theta,\rho-\rho_M)=n}\Gr^{M,\theta^\vee}.$$
The following facts are proved in \cite{BD, Section 5.3}.

\vskip3mm

\noindent{\bf Proposition.} {\it (a) The functor $\pi_!\gamma^*$ gives a map
$\res^{GM}\,:\,\Pc_G\to\tilde\Pc_M=\bigoplus_n\Pc_{M,n}[-n]$,
where $\Pc_{M,n}$ is the subcategory of $M(O)$-equivariant perverse sheaves
on $\Gr^{M,n}$.\hfill\break
(b) For any $\Ec,\Fc\in\Pc_G$ we have $\res^{GM}(\Ec\star\Fc)=(\res^{GM}\Ec)\star
(\res^{GM}\Fc)$.\hfill\break
(c) For any $\Ec\in\Pc_G$ we have $H^*(\Gr,\Ec)=H^*(\Gr^M,\res^{GM}\Ec)$.
\hfill\break
(d) If $P_1\subset P$ is a parabolic subgroup and $M_1$ is its Levi factor
then $\res^{MM_1}$ maps $\tilde\Pc_M$ to $\tilde\Pc_{M_1}$, and
$\res^{GM_1}=\res^{MM_1}\circ\res^{GM}$.}
\qed

\subhead 1.4\endsubhead
Let $\tilde\gen$ be the affine Kac-Moody Lie algebra associated to $G$.
Let $\tilde\omega_0$ be the fundamental weight of $\tilde\gen$ which is 
trivial on $\Lie(T)$. Let $W_0$ be the irreducible integrable highest weight
module of $\tilde\gen$ with higest weight $\tilde\omega_0$.
Let $\pi$ be the corresponding group homomorphism
$G(K)\to PGL(W_0)$ (see \cite{Ku, Appendix C} for instance).
The central extension $\tilde G(K)$ of $G(K)$ is the pull-back
$\pi^*GL(W_0)$, where $GL(W_0)$ must be viewed as a $\CC^\times$-principal
bundle on $PGL(W_0)$. The restriction of the central extension to $G(O)$,
denoted by $\tilde G(O)$, splits, i.e. $\tilde G(O)=G(O)\times\CC^\times$.
Fix a highest weight vector $w_0\in W_0$. 
Let $\Lc_G$ be the pull-back of $O_\PP(1)$ by the embedding of ind-schemes
$\iota\,:\,\Gr^G\hookrightarrow\PP(W_0)$
induced by the map
$$G(K)\to\PP(W_0),\, g\mapsto[\CC\cdot gw_0].$$
The sheaf $\Lc_G$ is obviously algebraic.

\subhead 1.5\endsubhead
For any $i\in I$ let $P_i$ be the corresponding subminimal parabolic
subgroup of $G$. Let $N_i\subset P_i$ be the unipotent radical and put 
$M_i=P_i/N_i$. Hereafter we set $i\res=\res^{GM_i}$, $\pi_i=\pi$, 
$\gamma_i=\gamma$, $Z_i=Z(M_i)$ and $\Lc_i=\Lc_{M_i}$. 
The product by the first Chern class of $\Lc_i$ gives a map 
$$l_i\,:\,H^*(\Gr^{M_i},\Ec)\to H^{*+2}(\Gr^{M_i},\Ec),$$
for any $\Ec\in\Pc_{M_i}$. 

\subhead 1.6\endsubhead
For any $\mu^\vee\in X^\vee$ set $S_{\mu^\vee}=U(K)\cdot e_{\mu^\vee}$.
It was proved by Mirkovic and Vilonen that if $\Ec\in\Pc_G$ then
$$H^*(\Gr,\Ec)=\bigoplus_{\mu^\vee\in X^\vee}
H_c^{2(\rho,\mu^\vee)}(S_{\mu^\vee},\Ec),\leqno(a)$$
(see \cite{MV}, and \cite{NP} for more details).
For any $i\in I$ and any $\mu^\vee\in X^\vee$ set also
$S^{M_i}_{\mu^\vee}=U_i(K)\cdot e_{\mu^\vee}\subset\Gr^{M_i}$.
The grassmanian $\Gr^{M_i}$ may be viewed as the set of points 
of $\Gr$ which are fixed by the action of the group $Z_i$ by left translations.
This fixpoints subset is denoted by ${}^{Z_i}\Gr$. In particular,
$S^{M_i}_{\mu^\vee}$ may be viewed as a subset of $\Gr$.

\head 2. Construction of the operators $\eb_i$, $\fb_i$, $\hb_i$\endhead
\subhead 2.1 \endsubhead
To avoid useless complications, hereafter we assume that $G$ is 
semi-simple. The generalization to the reductive case is immediate.
For any $i\in I$ and $\Ec\in\Pc_G$, let $\eb_i$ be the composition of the chain 
of maps
$$H^*(\Gr,\Ec)=H^*(\Gr^{M_i},i\res\Ec){\buildrel l_i\over\to}
H^{*+2}(\Gr^{M_i},i\res\Ec)=H^{*+2}(\Gr,\Ec).$$
Moreover, set
$$\hb_i=\bigoplus_{\lambda^\vee\in X^\vee}(\alpha_i,\lambda^\vee)\,
id_{H^*_c(S_{\lambda^\vee},\Ec)}\,:\,H^*(\Gr,\Ec)\to H^*(\Gr,\Ec).$$
By the hard Lefschetz theorem there is a unique linear operator 
$\fb_i\,:\,H^*(\Gr,\Ec)\to H^{*-2}(\Gr,\Ec)$ such that $(\eb_i,\hb_i,\fb_i)$
is a $\sen\len(2)$-triple.

\vskip3mm

\noindent{\bf Theorem.} {\it For any $\Ec\in\Pc_G$, the operators
$\eb_i,\fb_i,\hb_i$, with $i\in I$, give an action of the dual group $G^\vee$ on
the cohomology $H^*(\Gr,\Ec)$.}\qed

\subhead 2.2\endsubhead The rest of the paper is devoted to the proof of the 
theorem.

\vskip3mm

\noindent{\bf Lemma.} {\it For all $\lambda^\vee\in X^\vee$ we have
$$\eb_i(H^*_c(S_{\lambda^\vee},\Ec))\subset 
H^*_c(S_{\lambda^\vee+\alpha_i^\vee},\Ec)\and
\fb_i(H^*_c(S_{\lambda^\vee},\Ec))\subset 
H^*_c(S_{\lambda^\vee-\alpha_i^\vee},\Ec).$$
}

\vskip1mm

\noindent{\it Proof.}
It is sufficient to check the first claim.
Since 
$$S_{\lambda^\vee}=N_i(K)U_i(K)\cdot e_{\lambda^\vee}=
\pi_i^{-1}(S^{M_i}_{\lambda^\vee}),$$ 
we get, for any $\Ec\in\Pc_G$,
$$H^*_c(S_{\lambda^\vee},\Ec)=H^*_c(S^{M_i}_{\lambda^\vee},i\res\Ec).
\leqno(a)$$
Now, if $\Ec\in\Pc_{M_i}$ then
$$l_i(H^*_c(S^{M_i}_{\lambda^\vee},\Ec))=
l_i(H^{(\alpha_i,\lambda^\vee)}_c(S^{M_i}_{\lambda^\vee},\Ec))\subset
H^{2+(\alpha_i,\lambda^\vee)}_c(\Gr^{M_i},\Ec).$$
Moreover, for all $\mu^\vee\in X^\vee\simeq X_{T^\vee}$ we have
$$S_{\mu^\vee}\cap\Gr^{M_i,\theta^\vee}\neq\emptyset\,\iff\,
\mu^\vee|_{Z(M_i^\vee)}=\theta^\vee.$$
Thus, if $\Ec\in\Pc_{M_i}$ then
$$l_i(H^*_c(S^{M_i}_{\lambda^\vee},\Ec))=
\bigoplus_{\mu^\vee}H^{(\alpha_i,\mu^\vee)}_c(S^{M_i}_{\mu^\vee},\Ec),$$
where the sum is over all $\mu^\vee\in X^\vee\simeq X_{T^\vee}$ such that
$$\mu^\vee|_{Z(M_i^\vee)}=\lambda^\vee|_{Z(M_i^\vee)}\and
(\alpha_i,\mu^\vee)=(\alpha_i,\lambda^\vee+\alpha_i^\vee).$$
The only possibility is $\mu^\vee=\lambda^\vee+\alpha^\vee_i.$
\qed

\vskip3mm

\noindent The lemma implies that $[\hb_i,\eb_j]=(\alpha_i,\alpha_j^\vee)\,\eb_j$ 
for all $i,j\in I$. Since $\eb_i, \fb_i$, are locally nilpotent and since
$[\eb_i,\fb_i]=\hb_i$ by construction, if $[\eb_i,\fb_j]=0$ for any $i\neq j$ 
then the operators $\eb_i, \fb_i, \hb_i$,
give a representation of the Lie algebra $\gen^\vee$
of $G^\vee$ on the cohomology group $H^*(\Gr,\Ec)$ for any $\Ec\in\Pc_G$ 
(see \cite{Ka, Section 3.3}). 
The action of the operators $\hb_i$ lifts to an action of the torus of $G^\vee$.
Thus, the representation of the Lie algebra $\gen^\vee$ lifts to a representation of the group
$G^\vee$. By (1.3.d), in order to check the relation $[\eb_i,\fb_j]=0$ for $i\neq j$
we can assume that the group $G$ has rank 2.

\subhead 2.3\endsubhead Recall that any complex $\ICc_{\lambda^\vee}$ 
is a direct factor of a product
$\ICc_{\lambda_1^\vee}\star\ICc_{\lambda_2^\vee}\star\cdots\star
\ICc_{\lambda_n^\vee}$ such that the coweights $\lambda_i^\vee$ are 
either minuscule or quasi-minuscule (see \cite{NP, Proposition 9.6}). 
Observe that \cite{NP, Lemmes 10.2, 10.3} imply indeed that if the
set of minuscule coweights is non empty, then we can find such a product with all
the $\lambda^\vee_i$'s beeing minuscule.
Recall also that for any $\Ec,\Fc\in\Pc_G$ there is a canonical
isomorphism of graded vector spaces
$$H^*_c(S_{\lambda^\vee},\Ec\star\Fc)\simeq
\bigoplus_{\mu^\vee+\nu^\vee=\lambda^\vee}H^*_c(S_{\mu^\vee},\Ec)\otimes
H^*_c(S_{\nu^\vee},\Fc),\leqno(a)$$
(see \cite{MV}, and \cite{NP, Proof of Theorem 3.1} for more details).
Let $\Delta(\eb_i),\Delta(\fb_i),\Delta(\hb_i)$, be the composition
$$H^*(\Gr,\Ec)\otimes H^*(\Gr,\Fc)=H^*(\Gr,\Ec\star\Fc)
{\buildrel \eb_i,\fb_i,\hb_i\over\longrightarrow}
H^*(\Gr,\Ec\star\Fc)=H^*(\Gr,\Ec)\otimes H^*(\Gr,\Fc),$$
where the equalities are given by (1.6.a) and (a).

\vskip3mm

\noindent{\bf Lemma.} 
{\it If $x=\eb_i, \fb_i, \hb_i$, then $\Delta(x)=x\otimes 1+1\otimes x$.}

\vskip3mm

\noindent{\it Proof.} 
If $x=\hb_i$ the equation is obvious. 
If $x=\fb_i$ it is a direct consequence of the two others since a $\sen\len(2)$-triple
$(\eb_i,\hb_i,\fb_i)$ is completely determined by $\eb_i$ and $\hb_i$. 
Thus, from (2.3.a), (1.3.b) and (1.3.c), it suffices to check the 
equality when $G=\SL(2)$ and $x=\eb_i$. 
Then, the operator $\eb_i$ is the product by the 1-st Chern class of 
the line bundle $\Lc_{\SL(2)}$ on $\Gr^{\SL(2)}.$
More generally, for any simply connected group $G$, the $G(O)$-equivariant line bundle 
$\Lc_G$ on the grassmannian $\Gr$ lifts uniquely to a $G(K)$-equivariant line bundle
on the ind-scheme $G(K)\times_{G(O)}\Gr$. Let denote it by $\Lc_2$. 
The group $G(O)$ acts on the pull-back of $\Lc_G$ by the projection
$G(K)\times\Gr\to\Gr$. The quotient is the bundle $\Lc_2$. The vector bundle
$\Lc_2$ is algebraic, i.e. its restriction to
the subscheme 
$p^{-1}(\overline\Gr_{\lambda^\vee_1})\times_{G(O)}\overline\Gr_{\lambda_2^\vee}$ 
is an algebraic vector bundle
for any $\lambda_1^\vee$, $\lambda_2^\vee$.
Indeed, there is a normal pro-unipotent closed subgroup $H$ of $G(O)$
such that $G(O)/H$ is finite dimensional and $H$ acts trivialy on
$\overline\Gr_{\lambda^\vee_1},\overline\Gr_{\lambda^\vee_2}.$
Since $H$ is pro-unipotent, the restriction of $\Lc_G$ to 
$\overline\Gr_{\lambda^\vee_2}$ is $G(O)/H$-equivariant.
Thus the restriction of $\Lc_2$ to
$p^{-1}(\overline\Gr_{\lambda^\vee_1})\times_{G(O)}\overline\Gr_{\lambda_2^\vee}$ 
is identified with the algebraic sheaf on
$$(p^{-1}(\overline\Gr_{\lambda^\vee_1})/H)\times_{G(O)/H}\overline\Gr_{\lambda_2^\vee}$$ 
induced by the restriction of $\Lc_G$ to $\overline\Gr_{\lambda^\vee_2}$.
Consider also the pull-back
$\Lc_1$ of the line bundle $\Lc_G$ by the 1-st projection
$\tilde p\,:\,G(K)\times_{G(O)}\Gr\to\Gr$. We claim that 
$$m^*\Lc_G=\Lc_1\otimes\Lc_2.\leqno(b)$$
Let $\mu\,:\,G(K)\times G(K)\to G(K)$ be the multiplication map. 
The product in the group $\tilde G(K)$ gives an isomorphism of bundles
$$\mu^*p^*\Lc_G\simeq p^*\Lc_G\boxtimes p^*\Lc_G$$
on $G(K)\times G(K)$. This isomorphism descends to the fiber product
$G(K)\times_{G(O)}\Gr$ and implies (b). 
Observe now that (a) is induced by the canonical isomorphism
$$\bigl(p^{-1}(\overline\Gr_{\lambda^\vee_1})\times_{G(O)}
\overline\Gr_{\lambda^\vee_2}\bigr)\cap m^{-1}(S_{\lambda^\vee})\simeq
\bigsqcup_{\mu^\vee+\nu^\vee=\lambda^\vee}
(S_{\mu^\vee}\cap\overline\Gr_{\lambda_1^\vee})\times
(S_{\nu^\vee}\cap\overline\Gr_{\lambda_2^\vee})$$
resulting from the local triviality of $p$. Let $l_\Ec$
be the product by the 1-st Chern class of $\Lc_G$ on the global
cohomology of the perverse sheaf $\Ec\in\Pc_G$. 
Then (a) and (b) give $l_{\Ec\star\Fc}=l_\Ec\otimes 1+1\otimes l_\Fc$.
\qed

\subhead 2.4\endsubhead From Section 2.3 and (1.3.d)
we are reduced to check the 
relation $[\eb_i,\fb_j]=0$, $i\neq j$, on the cohomology group $H^*(\Gr,\Ec)$
when $G$ has rank 2 and $\Ec=\ICc_{\lambda^\vee}$, with $\lambda^\vee$
minuscule or quasi-minuscule. For any dominant coroot $\lambda^\vee$
let $\Omega(\lambda^\vee)\subset X^\vee$ be the set of weights
of the simple $G^\vee$-module with highest weight $\lambda^\vee$. Recall that

\itemitem{(a)} the coweight $\lambda^\vee\in X^\vee_+-\{0\}$ is minuscule 
if and only if $\Omega(\lambda^\vee)=W\cdot\lambda^\vee$, if and only if
$(\alpha,\lambda^\vee)=0,\pm 1$, for all $\alpha\in R$,

\itemitem{(b)} the coweight $\lambda^\vee\in X^\vee_+-\{0\}$ is 
quasi-minuscule if and only
if $\Omega(\lambda^\vee)=W\cdot\lambda^\vee\cup\{0\}$, if and only if
$\lambda^\vee$ is a maximal short coroot. Moreover if $\lambda^\vee$ is 
quasi-minuscule then $(\alpha,\lambda^\vee)=0,\pm 1,$ for all 
$\alpha\in R-\{\pm\lambda\}$.

\noindent 
For any coweight $\lambda^\vee$ we consider the 
isotropy subgroup $G_{\lambda^\vee}$ of $e_{\lambda^\vee}$ in $G$. 
Thus
$$G_{\lambda^\vee}=T\prod_{(\alpha,\lambda^\vee)\leq 0}U_\alpha.$$ 
In particular $B^-\subset G_{\lambda^\vee}$ if $\lambda^\vee$ is dominant,
and we can consider the line bundle $\Lc(\lambda)$ on $G/G_{\lambda^\vee}$ 
associated to the weight $\lambda$. The structure of $\overline\Gr_{\lambda^\vee}$
for $\lambda^\vee$ minuscule or quasi-minuscule is described as follows in \cite{NP}.

\vskip3mm

\noindent{\bf Proposition.} {\it
(c) If $S_{\mu^\vee}\cap\overline\Gr_{\lambda^\vee}\neq\emptyset$, then $\mu^\vee\in\Omega(\lambda^\vee)$.
\hfill\break
(d) If $\mu^\vee\in W\cdot\lambda^\vee$, then 
$S_{\mu^\vee}\cap\overline\Gr_{\lambda^\vee}=S_{\mu^\vee}\cap\Gr_{\lambda^\vee}$.
\hfill\break
(e) If $\lambda^\vee\in\X_+^\vee$ is minuscule, then
$$\overline\Gr_{\lambda^\vee}=\Gr_{\lambda^\vee}=G/G_{\lambda^\vee}
\and S_{w\cdot\lambda^\vee}\cap\Gr_{\lambda^\vee}\simeq
UwG_{\lambda^\vee}/G_{\lambda^\vee}\quad\forall w\in W.$$
(f) Assume that $\lambda^\vee\in X^\vee_+$ is quasi-minuscule. 
Then $\Gr_{\lambda^\vee}\simeq\Lc(\lambda)$ and
$\overline\Gr_{\lambda^\vee}\simeq\Lc(\lambda)\cup\{e_0\}$ 
as a $G$-varieties. Moreover,
$$S_{w\cdot\lambda^\vee}\cap\Gr_{\lambda^\vee}\simeq\left\{
\matrix
UwG_{\lambda^\vee}/G_{\lambda^\vee}\quad\hfill&\text{if}\quad w\cdot\lambda\in R_-,\hfill\cr\cr
\Lc|_{UwG_{\lambda^\vee}/G_{\lambda^\vee}}\quad\hfill&\text{if}\quad w\cdot\lambda\in R_+.\hfill
\endmatrix\right.$$
}\qed

\head 3. Proof of the relation $[\eb_i,\fb_j]=0$\endhead
\subhead 3.1\endsubhead Assume that $G$ has rank two and set $I=\{1,2\}$.
The Bruhat decomposition for $M_i$ implies that 
$M_ie_{\lambda^\vee}=U_ie_{s_i\cdot\lambda^\vee}
\cup U_iU_{-\alpha_i^\vee}e_{\lambda^\vee}.$
Thus,

\itemitem{(a)} if $(\alpha_i,\lambda^\vee)>0$ then
$M_ie_{\lambda^\vee}=M_ie_{s_i\cdot\lambda^\vee}=
U_ie_{\lambda^\vee}\cup\{e_{s_i\cdot\lambda^\vee}\},$

\itemitem{(b)} if $(\alpha_i,\lambda^\vee)=0$ then
$M_ie_{\lambda^\vee}=\{e_{\lambda^\vee}\}$.

\subhead 3.2\endsubhead 
Assume that $\lambda^\vee$ is a minuscule dominant coweight.
Fix $\mu^\vee=w\cdot\lambda^\vee$ with $w\in W$, 
and fix $i\in I$. One of the following three cases holds 

\itemitem{(a)} we have $(\alpha_i,\mu^\vee)=1$, and 
$$S_{\mu^\vee}^{M_i}\cap\Gr_{\lambda^\vee}=U_ie_{\mu^\vee},\quad
S_{\mu^\vee-\alpha_i^\vee}^{M_i}\cap\Gr_{\lambda^\vee}=
\{e_{\mu^\vee-\alpha_i^\vee}\},\quad\Gr^{M_i}_{\mu^\vee}=
U_ie_{\mu^\vee}\cup\{e_{\mu^\vee-\alpha_i^\vee}\},$$

\itemitem{(b)} we have $(\alpha_i,\mu^\vee)=-1$, and
$$S_{\mu^\vee+\alpha_i^\vee}^{M_i}\cap\Gr_{\lambda^\vee}=
U_ie_{\mu^\vee+\alpha_i^\vee},\quad
S_{\mu^\vee}^{M_i}\cap\Gr_{\lambda^\vee}=\{e_{\mu^\vee}\},\quad
\Gr^{M_i}_{\mu^\vee}=U_ie_{\mu^\vee+\alpha_i^\vee}\cup 
\{e_{\mu^\vee}\},$$

\itemitem{(c)} we have $(\alpha_i,\mu^\vee)=0$, and 
$$S_{\mu^\vee}^{M_i}\cap\Gr_{\lambda^\vee}=
\Gr^{M_i}_{\mu^\vee}=\{e_{\mu^\vee}\}.$$

\vskip2mm

\noindent 
Obviously, the sheaf $i\res\Ic\Cc_{\lambda^\vee}$ is supported on
$\Gr_{\lambda^\vee}\cap\Gr^{M_i}$. Thus, by (2.2.a) and Lemma 2.2, if 
$\fb_2\eb_1(H^*_c(S_{\mu^\vee},\ICc_{\lambda^\vee}))\neq\{0\}$ then
$$S_{\mu^\vee}^{M_1}\cap\Gr_{\lambda^\vee},\quad
S_{\mu^\vee+\alpha_1^\vee}^{M_1}\cap\Gr_{\lambda^\vee},\quad
S_{\mu^\vee+\alpha_1^\vee}^{M_2}\cap\Gr_{\lambda^\vee},\quad
S_{\mu^\vee+\alpha_1^\vee-\alpha_2^\vee}^{M_2}\cap\Gr_{\lambda^\vee}$$
are non empty.
In particular, we get 
$$(\alpha_1,\mu^\vee)=-1,\qquad (\alpha_2,\mu^\vee+\alpha_1^\vee)=1.$$
Since $\lambda^\vee$ is minuscule, $\mu^\vee\in W\cdot\lambda^\vee$, and $(\alpha_2,\alpha^\vee_1)\leq 0$,
we get 
$$(\alpha_2,\mu^\vee)=1,\quad (\alpha_1,\mu^\vee)=-1,\quad (\alpha_2,\alpha^\vee_1)=0.$$
Similarly, 
$$\eb_1\fb_2(H^*_c(S_{\mu^\vee},\ICc_{\lambda^\vee}))\neq\{0\}\Rightarrow
(\alpha_2,\mu^\vee)=1,\quad (\alpha_1,\mu^\vee)=-1,\quad
(\alpha_1,\alpha^\vee_2)=0.$$
Thus we are reduced to the case where 
$G=\SL(2)\times\SL(2)$, $M_1\simeq\SL(2)\times\{1\}$,
$M_2=\{1\}\times\SL(2)$, $\lambda^\vee=\omega_1^\vee+\omega_2^\vee$,
$\mu^\vee=-\omega_1^\vee+\omega_2^\vee$, and $\Ic\Cc_{\lambda^\vee}$
is the constant sheaf on $\Gr_{\lambda^\vee}$. Then,
$$\Gr_{\lambda^\vee}\simeq\PP^1\times\PP^1,\quad
\Gr^{M_1}\cap\Gr_{\lambda^\vee}\simeq\PP^1\times\{0,\infty\},\quad
\Gr^{M_2}\cap\Gr_{\lambda^\vee}\simeq\{0,\infty\}\times\PP^1.$$
Recall that, with the notations of Section 1.4, the fiber of $\Lc_G^{-1}$ at 
$e_{\lambda^\vee}$ is identified with $\CC t^{\lambda^\vee}w_0$.
Recall also that the extended affine Weyl group $W\ltimes X^\vee$ 
acts on the lattice $\Hom(T\times\GG_m,\GG_m)$ in such a way that 
$\lambda^\vee\cdot\tilde\omega_0=\lambda+\tilde\omega_0$ 
for all $\lambda^\vee\in X^\vee$ 
(see \cite{PS, Proposition 4.9.5} for instance).
Thus, for any dominant coweight
$\lambda^\vee$ the restriction of $\Lc_G$ to the $G$-orbit 
$G\cdot e_{\lambda^\vee}$ is the line bundle 
$\Lc(\lambda)$ on $G/G_{\lambda^\vee}$.
In particular the restriction of the line bundle $\Lc_i$ to 
$\Gr^{M_i}_{\lambda^\vee}$ is $O_{\PP^1}(1)$.
Thus $\eb_1=l\otimes id$ and $\eb_2=id\otimes l$,
where $l$ is the product by the 1-st Chern class of $O_{\PP^1}(1)$.
The relation is obviously satisfied.

\subhead 3.3\endsubhead Assume that $\lambda^\vee$ is a
quasi-minuscule dominant coweight.
Observe that if $G$ is of type $A_1\times A_1$, $A_2$ or $B_2$,
then the set of minuscule coweights is non empty.
Thus, from Section 2.3 we can assume that $G$ is of type $G_2$.
Let $\alpha_1^\vee$ be the long simple coroot, and let $\alpha_2^\vee$ be the short one.
Then
$$\lambda^\vee=\alpha_1^\vee+2\alpha_2^\vee,\qquad(\alpha_2,\alpha_1^\vee)=-3,\qquad
(\alpha_1,\alpha_2^\vee)=-1.$$
Set $\Lc=\Lc(\lambda)$ and $\bar{\Lc}=\Lc\cup\{e_0\}$. 
Then 
$\overline\Gr_{\lambda^\vee}\cap\Gr^{M_i}$ is the fixpoints set of $Z_i$ 
on $\bar{\Lc}$, i.e.
$$\overline\Gr_{\lambda^\vee}\cap\Gr^{M_i}=\{e_0\}\cup
\bigcup_{\mu^\vee\in W\cdot\lambda^\vee}\Gr^{M_i}_{\mu^\vee}
\quad\text{where}\quad
\Gr^{M_i}_{\mu^\vee}={}^{Z_i}\Lc|_{M_ie_{\mu^\vee}}.$$
Assume that $\mu^\vee=w\cdot\lambda^\vee$ with $w\in W$. 

\itemitem{(a)} If $(\alpha_i,\mu^\vee)=0$ then 
$\Gr^{M_i}_{\mu^\vee}={}^{Z_i}\Lc|_{e_{\mu^\vee}}$.
The torus $T$ acts on the fiber $\Lc|_{e_{\mu^\vee}}$ by the character $\mu$.
Since $\mu^\vee\neq 0$ and $(\alpha_i,\mu^\vee)=0$, necessarily
$\mu(Z_i)$ is non-trivial. Thus, $\Gr^{M_i}_{\mu^\vee}=e_{\mu^\vee}.$

\itemitem{(b)} If $(\alpha_i,\mu^\vee)=2$ then $\mu^\vee=\alpha_i^\vee$
and $\Gr^{M_i}_{\mu^\vee}=\Lc|_{M_ie_{\alpha_i^\vee}}.$
Moreover, since $\lambda^\vee$ is short and $\alpha_1^\vee$ is long
we have $i=2$.

\itemitem{(c)} If $(\alpha_i,\mu^\vee)=1$ then 
$\Gr^{M_i}_{\mu^\vee}=M_ie_{\mu^\vee}$ because 
$$\bigl(\mu(Z_i)=1\and\mu^\vee\in R^\vee\bigr)\,\Rightarrow
\mu^\vee\in\ZZ\alpha_i^\vee\,\Rightarrow(\alpha_i,\mu^\vee)\neq 1.$$

\noindent 
In Case (b) we get ($i=2$) 
$$S^{M_i}_{\alpha_i^\vee}\cap\overline\Gr_{\lambda^\vee}=
\Lc|_{U_ie_{\alpha_i^\vee}},\quad
S^{M_i}_{-\alpha_i^\vee}\cap\overline\Gr_{\lambda^\vee}=
e_{-\alpha_i^\vee},\quad
S^{M_i}_0\cap\overline\Gr_{\lambda^\vee}=
\bar{\Lc}^\times|_{e_{-\alpha_i^\vee}},$$
$$\and\overline\Gr_{\alpha_i^\vee}^{M_i}=\bar{\Lc}|_{M_ie_{\alpha_i^\vee}},$$
where the upperscript $\times$ means than the zero section has been removed.
In Case (c) we get 
$$S_{\mu^\vee}^{M_i}\cap\overline\Gr_{\lambda^\vee}=U_ie_{\mu^\vee},\quad
S_{\mu^\vee-\alpha_i^\vee}^{M_i}\cap\overline\Gr_{\lambda^\vee}=
e_{\mu^\vee-\alpha_i^\vee},\and\overline\Gr_{\mu^\vee}^{M_i}=
U_ie_{\mu^\vee}\cup e_{\mu^\vee-\alpha_i^\vee}.$$

\vskip2mm

\noindent Thus, for any $\mu^\vee\in X^\vee$, Claim
(2.2.a) and Lemma 2.2 imply that 

\vskip2mm

\noindent(d) if $\eb_1(H^*_c(S_{\mu^\vee},\ICc_{\lambda^\vee}))\neq\{0\}$ then
$(\alpha_1,\mu^\vee)=-1$, or $\mu^\vee=0$, or $\mu^\vee=-\alpha_1^\vee$,

\vskip2mm

\noindent(e) if $\fb_2(H^*_c(S_{\mu^\vee},\ICc_{\lambda^\vee}))\neq\{0\}$ then
$(\alpha_2,\mu^\vee)=1$, or $\mu^\vee=0$, or $\mu^\vee=\alpha_2^\vee$.

\vskip2mm

\noindent 
Observe that in Case (d) the identity (2.4.c) and Lemma 2.2 imply indeed that $\mu^\vee\neq 0,-\alpha_1^\vee$,
because
$$H^*_c(S_{\alpha_1^\vee},\ICc_{\lambda^\vee})=H^*_c(S_{-\alpha_1^\vee},\ICc_{\lambda^\vee})
=\{0\}.$$
Thus, if $\fb_2\eb_1(H^*_c(S_{\mu^\vee},\ICc_{\lambda^\vee}))\neq\{0\}$ then
$(\alpha_1,\mu^\vee)=-1$ and $(\alpha_2,\mu^\vee+\alpha_1^\vee)=1.$
We get $(\alpha_2,\mu^\vee)=4$. This is not possible since $\mu^\vee\in\Omega(\lambda^\vee)$
and $\lambda^\vee$ is quasi-minuscule.
Similarly, if $\eb_1\fb_2(H^*_c(S_{\mu^\vee},\ICc_{\lambda^\vee}))\neq 0$ then
$(\alpha_1,\mu^\vee)=-2$. This is not possible either.
Thus, the relation $[\eb_1,\fb_2]=0$ is obviously satisfied.
The relation $[\eb_2,\fb_1]=0$ is proved in the same way.

\vskip1cm

\noindent{\it Acnowledgements.}
{\eightpoint{The author would like to thank A. Beilinson and B.C. Ng\^o for discussion on this subject.}}

\vskip3cm
{\eightpoint{
$$\matrix\format&\l&\l&\l\\
\phantom{.} & {\text{Eric Vasserot}}\\
\phantom{.} & {\text{D\'epartement de Math\'ematiques}}\\
\phantom{.} & {\text{Universit\'e de Cergy-Pontoise}}\\
\phantom{.} & {\text{2 Av. A. Chauvin}}\\
\phantom{.} & {\text{95302 Cergy-Pontoise Cedex}}\\
\phantom{.} & {\roman{France}}\\
            & {\text{email: eric.vasserot\@u-cergy.fr}}
\endmatrix$$
}}
\enddocument